\documentclass{notices}
\newcommand{\vect}[1]{\boldsymbol{{#1}}}
\usepackage{amsfonts}
\usepackage{amssymb}
\usepackage{graphicx,color}

\title{Applications of PDEs and Stochastic Modeling to Protein Transport in Cell Biology 
}

\author{
  Maria-Veronica Ciocanel
  \affil{
    ciocanel@math.duke.edu.
    }
}

\begin{document}

\maketitle

\section*{Introduction}

Intracellular transport refers to the movement of proteins and vesicles inside cells. These transport processes are essential to the healthy development of many organisms as well as more generally to healthy cellular function. Inside cells, various proteins and protein filaments must be dynamic and interact on robust time and spatial scales to ensure proper cell functioning. In particular, healthy sorting and delivery of protein components relies on the interaction of several molecular players: the cargo proteins (which need to be carried and delivered), the molecular motors (the protein vehicles that actively move the cargo), and the cytoskeleton network (the protein filaments that act as roads for cargo transport). 

There are many aspects that complicate intracellular transport dynamics. One aspect is that there are several types of cytoskeleton filaments that are needed for different cell requirements. Microtubules interact with motor proteins kinesin and dynein, give the cell structure, and provide long highways for transport of proteins. Actin filaments interact with myosin motor proteins and are well-known for their role in muscle contraction, but they also contribute to stability and motion in most cells. In addition, molecular motors move in different directions along these filaments (for example, kinesin and dynein move in opposite directions along microtubules), thus leading to bidirectional transport of proteins along these filament networks. These motors have also been found to be co-dependent, thus leading to interesting competing hypotheses on their interactions \cite{hancock2014bidirectional}. The complex dynamics and interactions between protein molecules and filaments on different time and spatial scales generate many opportunities for mathematical modeling and analysis that can provide insights into protein sorting and overall cellular organization. 

In addition, experimental techniques in cell biology have undergone great advances in recent years. For example, fluorescence microscopy and live imaging allow researchers to observe the localization of various proteins inside cells, while optical trapping experiments can apply mechanical forces on molecules and thus provide insights into molecular motors and their motion. However, it is difficult to observe all key protein players and their interactions using experiments, especially in living systems, where these observations would be of utmost physiological relevance. For example, in many \textit{in vivo} studies, protein cargo trajectories may be observed, but the underlying cytoskeleton network locations are unknown. These limitations once again provide exciting possibilities for mathematical modeling and analysis of the interplay between protein dynamics and filament geometry, which in turn can inform additional experimental directions.

\section*{PDE Modeling of Intracellular Transport}

In intracellular transport, protein cargo may need to be localized to a specific domain of the cell, or may have to be uniformly spread throughout cellular processes. For example, messenger RNA (mRNA) must localize at the vegetal (bottom) side of the developing egg cell in the frog \textit{Xenopus laevis}, while neurofilaments (intermediate protein polymers) are typically dynamic and spread out throughout the axons of neurons in mice and rats. The timing of these transport dynamics is also key: mRNA localization in the developing frog consistently takes 24-48 hours to achieve. 

Understanding spatiotemporal cargo dynamics is therefore important for addressing mechanistic questions inside cells, making partial differential equations (PDEs) a useful modeling tool. This type of modeling requires a choice of the number of spatial dimensions to account for. In systems such as neuronal axons or fungal hyphae, transport is often assumed to be one-dimensional given the geometry of the cells and the parallel orientation of the cytoskeletal filaments. On the other hand, systems such as developing oocytes or the budding yeast are often represented in two dimensions for simplicity, and occasionally in three dimensions if there is available data to support that level of complexity \cite{trong2015cortical}.

PDE models have been used to model aspects of intracellular transport for a few decades. Early work proposed systems of PDEs to model the mechanical cycle of kinesin \cite{peskin1995coordinated} and the fast axonal transport of organelles, whose dynamics was characterized using approximate traveling waves \cites{blum1985model,reed1990approximate}. Studies \cites{brooks1999probabilistic,friedman2006approximate,friedman2007uniform} then rigorously validated this approximation for a range of PDE models of active transport. Other directions focused on models of molecular cargo driven by different numbers of motor proteins or multiple motor families \cites{klumpp2005cooperative,muller2008tug,kunwar2010robust,muller2010bidirectional,arpaug2019motor}. Motivated by the delivery of cargo to specific sites (such as a target synapse), others interpreted the motor-driven transport dynamics as an intermittent trapping process \cite{newby2010local,newby2010quasi}. PDE modeling frameworks have led to insights into microtubule and pigment granule self-organization \cite{cytrynbaum2004computational}, into stop-and-go motion of neurofilament proteins in axons \cites{jung2009modeling,li2012axonal}, and into regulation of early endosome transport in fungal cells \cite{gou2014mathematical,dauvergne2015application}, to name a few. The review \cite{mogre2020getting} provides a detailed overview of the biological objectives of intracellular transport and of models describing diffusion, motor-driven transport, and cytoplasmic flow mechanisms for accomplishing intracellular transport.

\begin{figure*}
\centering
\includegraphics[width=1\linewidth]{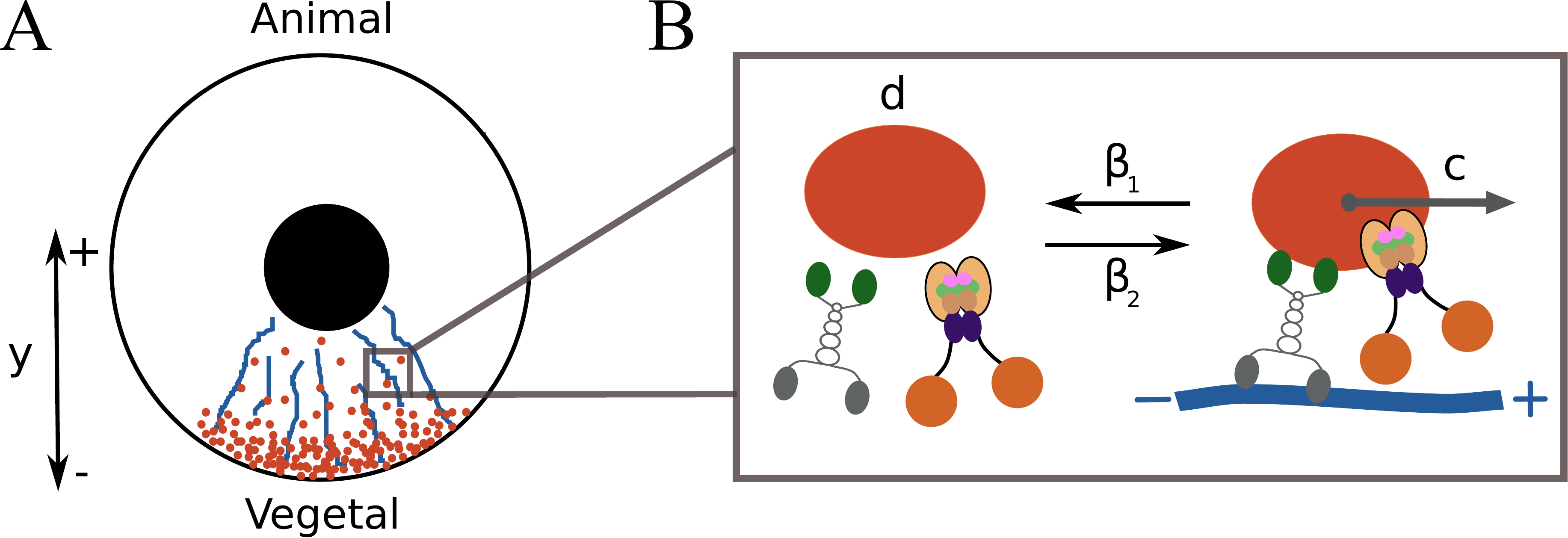}
\caption{(A) Illustration of an example application in intracellular transport: the localization of mRNA molecules (red particles) along microtubule filaments (blue lines) during development of the frog \textit{Xenopus laevis}. This establishes an animal-vegetal axis of development for the organism. (B) Simple 2-state model describing switching between diffusion of mRNAs (with diffusion coefficient $d$) and active transport of the particles driven by molecular motors along a microtubule filament (with speed $c$).}
\label{fig:mRNA_oocyte_2state}
\end{figure*}

In what follows, we will consider the motivating example of mRNA transport in developing frog oocytes, illustrated in Figure~\ref{fig:mRNA_oocyte_2state}A. In many cell types, mRNA molecules accumulate to specific regions of the cell and are thus asymmetrically distributed. This asymmetric segregation of maternal mRNA molecules allows for select proteins to be produced in specific cell locations, and ultimately produces a top-down (animal-vegetal) polarity of the cell, which establishes the body plan for the embryo and is crucial for normal development. 

Diffusion of mRNAs (as well as other protein cargos) has been observed in various organisms, and experimental techniques such as fluorescence recovery after photobleaching (FRAP) can be used to determine whether diffusion contributes to the dynamics \cite{sprague2004analysis}. Experiments knocking down microtubules or interfering with molecular motors \textit{in vivo} have also shown that active transport along microtubules is required for mRNA localization \cite{gagnon2013directional,messitt2008multiple}. mRNA cargo may also be driven by several motors with different speeds (and direction preferences) along microtubules. These proteins have also been observed to pause in certain cell measurements, due to various potential underlying mechanisms \cite{hancock2014bidirectional}. We therefore describe the dynamics of the system using systems of advection-reaction-diffusion partial differential equations incorporating all known biological dynamics, as well as plausible switching between states. 

For instance, a simple model describing mRNA molecules that transition between active movement and diffusion is illustrated in Figure~\ref{fig:mRNA_oocyte_2state}B and can be described using the system:
\begin{align}\label{eq:2-state}
    \partial_t{u} &= c \partial_y{u} - \beta_1 u + \beta_2 v \nonumber \\
    \partial_t{v} &= d \nabla^2 v + \beta_1 u - \beta_2 v\,,
\end{align}
with $\nabla = \begin{pmatrix} \partial_x \\ \partial_y \end{pmatrix}$ and $\nabla^2 = \partial^2_x + \partial^2_y $. Here $y$ denotes the spatial dimension from the nucleus to the cell cortex \cite{MR3892411}, as indicated in Figure~\ref{fig:mRNA_oocyte_2state}A. Variable $u(x,y,t)$ corresponds to the concentration of mRNA that moves to the cortex with speed $c$ (presumably driven by molecular motors along a microtubule) and $v(x,y,t)$ corresponds to the concentration of mRNA that diffuses in the cell cytoplasm with diffusion coefficient $d$. Parameter $\beta_1$ denotes the rate of transition from the moving to the diffusing state, while $\beta_2$ denotes the rate of transition from the diffusing to the moving state. Here, advection is assumed to occur only in the $y$ direction, while diffusion is 2-dimensional.

In general, protein cargo transport in different biological systems may be better characterized by models with additional dynamic states. Equations~\eqref{eq:2-state} can therefore be generalized to:
\begin{align}\label{eq:gen}
\partial_t\begin{pmatrix} u_1\\ ... \\ u_n  \end{pmatrix} = D \nabla^2\begin{pmatrix}  u_1\\ ... \\ u_n  \end{pmatrix} + C \partial_y\begin{pmatrix} u_1\\ ... \\ u_n  \end{pmatrix} + A \begin{pmatrix} u_1\\ ... \\ u_n \end{pmatrix}
\end{align}
where each variable $u_j$ models the concentration of particles in a certain state. The speeds $c_j$ and diffusion coefficients $d_j$ for each state are stored in diagonal matrices $C = \mathrm{diag}(c_j)$ and $D=\mathrm{diag}(d_j)$, while the transitions between states are described by the reaction-rate matrix $A$. For the simple 2-state model in equations~\eqref{eq:2-state}, the parameter matrices consist of:
\begin{align*}
    D &= \begin{pmatrix} 0 & 0 \\ 0 & d \end{pmatrix} \,,\\
    C &= \begin{pmatrix} c & 0 \\ 0 & 0 \end{pmatrix} \,,\\
    A &= \begin{pmatrix} -\beta_1 & \beta_2 \\ \beta_1 & -\beta_2 \end{pmatrix} \,.
\end{align*}

While equations such as \eqref{eq:gen} are not analytically tractable, we can use tools from dynamical systems to understand the protein dynamics in the asymptotic limit of large time. This limit is appropriate since transport processes take place on relatively large developmental timescales, such as 1-2 days for mRNA localization. For simplicity, we outline the analysis for the case of one-dimensional diffusion in the $y$ direction, as in \cite{ciocanel2017analysis}. We consider the Fourier transform ansatz 
\begin{align} \label{eq:ansatz}
(u_1, ...,  u_n)^T &= e^{\lambda t + iky} \vect{\tilde{u}_\mathrm{ic}}\,,
\end{align}
where $\vect{\tilde{u}_\mathrm{ic}}$ is a vector of initial conditions and $k$ is the wavenumber. Setting $\nu=ik$ and evaluating this ansatz in equation~\eqref{eq:gen} leads to
\begin{align} \label{eq:lambda_nu}
(A + \nu C + \nu^2D - \lambda I) \vect{\tilde{u}_\mathrm{ic}} = 0\,,
\end{align}
which is equivalent to $L(\lambda,\nu)\vect{\tilde{u}_\mathrm{ic}}=0$ for operator $L(\lambda,\nu) = A + \nu C + \nu^2D - \lambda I$, with $L(0,0)=A$. 
Since $A$ models conservation of particles in these models, we can verify that operator $L$ satisfies the conditions for applying Lyapunov-Schmidt reduction theory to study the solutions of \eqref{eq:lambda_nu} \cite{ciocanel2017analysis,MR3892411}. By projecting the equation onto the range and kernel of the operator, we find that
\begin{align}\label{eq:lamb_expansion}
    \lambda &= \nu \frac{\langle \vect{\psi_0},C\vect{u_0}\rangle}{\langle \vect{\psi_0}, \vect{u_0} \rangle} + \nu^2 \frac{\langle \vect{\psi_0},(D-C\tilde{A}^{-1}\tilde{C})\vect{u_0}\rangle}{\langle \vect{\psi_0}, \vect{u_0} \rangle}\,,
\end{align}
where $\vect{u_0}$ is the eigenvector of the zero eigenvalue of reaction matrix $A$, $\vect{\psi_0}$ is the eigenvector of the zero eigenvalue of the adjoint matrix $A^*$, $\tilde{A}$ is the restriction of transition matrix $A$ to its range, $\tilde{C}=C-\frac{\langle \vect{\psi_0},C\vect{u_0}\rangle}{\langle \vect{\psi_0}, \vect{u_0} \rangle}I$, and $\langle , \rangle$ denotes the dot product. 

We return to the ansatz in \eqref{eq:ansatz} and assume a Dirac delta function initial condition $u_{\mathrm{ic},l} = \delta(y)$, modeling a single particle starting at $y=0$. Then we take the inverse Fourier transform and the asymptotic limit of large time $t$ to obtain the spreading Gaussian solution
\begin{align}\label{eq:ul}
   u_l(y,t) = \frac{1}{\sqrt{2\pi \sigma_{\mathrm{eff}} t}} e^{-\frac{(y+v_{\mathrm{eff}}t)^2}{2\sigma_{\mathrm{eff}} t}}
\end{align}
for each cargo population $u_l$. The solution is therefore characterized by the effective velocity $v_{\mathrm{eff}}$ and effective diffusivity $\sigma_{\mathrm{eff}}$, which are given by:
\begin{align}\label{eq:effective_qtities}
    v_{\mathrm{eff}} &= \frac{\langle \vect{\psi_0},C\vect{u_0}\rangle}{\langle \vect{\psi_0}, \vect{u_0} \rangle}  \,, \nonumber \\
    \sigma_{\mathrm{eff}} &= \frac{\langle \vect{\psi_0},(D-C\tilde{A}^{-1}\tilde{C})\vect{u_0}\rangle}{\langle \vect{\psi_0}, \vect{u_0} \rangle}\,.
\end{align}
These quantities describe the evolution of the total mRNA population (i.e., $\sum_{l=1}^n u_l$) and its dependence on model parameters. This is especially important since most fluorescence microscopy experiments cannot distinguish between mRNA in different dynamical states, so that describing the effective transport of the proteins is very useful. For the simple 2-state model in \eqref{eq:2-state}, we derive the effective velocity and diffusion of the protein cargo as:
\begin{align}\label{eq:effective_qtities_2state}
    v_{\mathrm{eff}} &= c \frac{\beta_2}{\beta_1+\beta_2}  \,, \nonumber \\
    \sigma_{\mathrm{eff}} &= d \frac{\beta_1}{\beta_1+\beta_2} + c^2 \frac{\beta_1\beta_2}{(\beta_1+\beta_2)^3}  \,.
\end{align}

This dynamical systems approach to the asymptotic analysis of equations of the form \eqref{eq:gen} matches results of prior analysis of Chapman-Kolmogorov equations for the probability densities of particles engaged in a tug-of-war model of motor-driven cargo transport \cite{newby2010random}. Similar PDE systems that describe the evolution of the probabilities of a particle to be in different velocity or diffusive states have been reviewed in \cite{bressloff2013stochastic}. Previous work has reduced such PDE systems to a scalar Fokker-Planck equation using quasi-steady-state methods, which assume that the state transition rates are faster than the velocity of moving states on relevant timescales. As discussed in the next section, these parameter assumptions may not always hold when studying intracellular transport dynamics. The approach described here therefore complements prior work by considering a large time limit with no additional parameter assumptions.

\subsection*{Connecting models with fluorescence microscopy data}\label{sec:fluorescence_data}

Determining appropriate parameters for dynamical systems models is a common challenge in many mathematical biology applications. When it comes to questions about transport inside cells, biologists are very interested in estimates for parameter values in models such as \eqref{eq:gen}. Some estimates of protein diffusivities or active speeds may be informed from \textit{in vitro} experiments, but the crowded living cellular environment can significantly modify this dynamics. Similarly, the dynamics of protein complexes (such as motor-cargo complexes) can be more difficult to estimate, and switching rates between dynamic states are often unknown. 

\begin{figure*}
\centering
\includegraphics[width=01\linewidth]{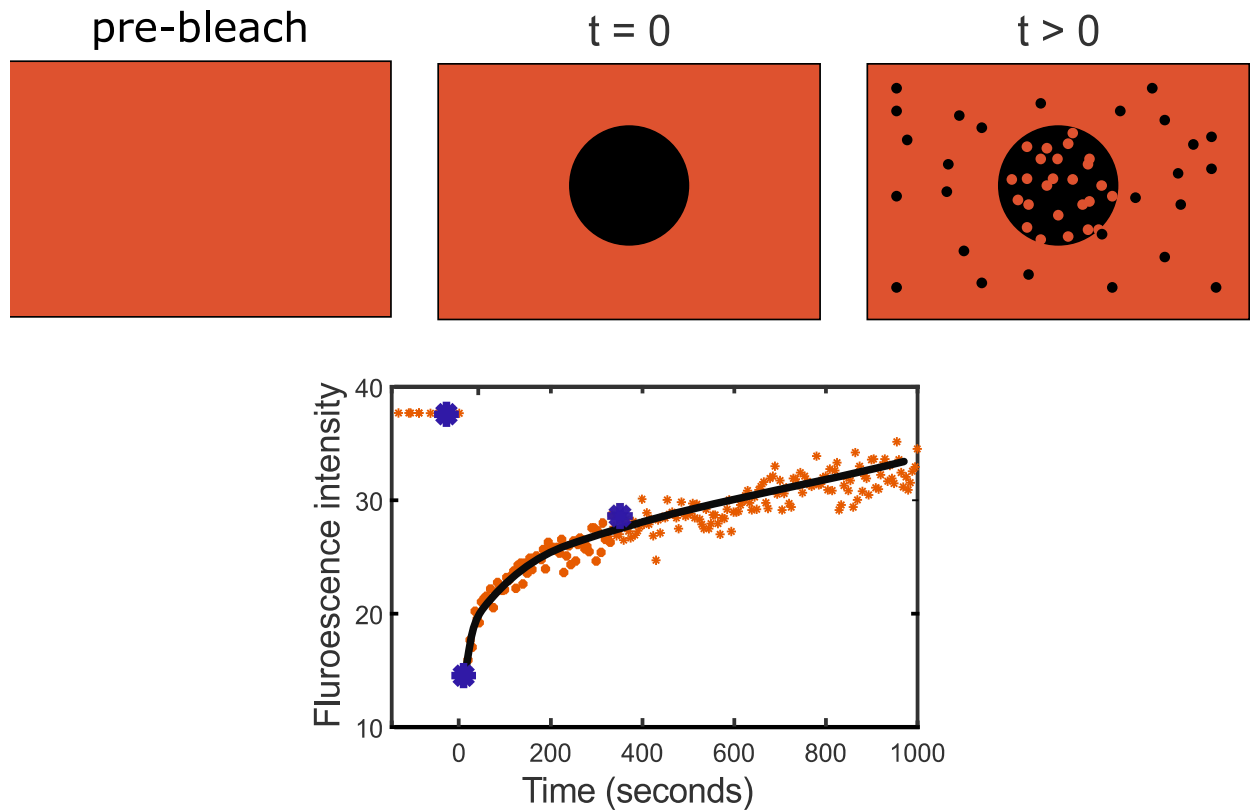}
\caption{Top: Fluorescence Recovery After Photobleaching (FRAP) involves bleaching of a small spot in a fluoresced region of the cell. At later times, non-bleached and fluorescent molecules mix between the regions. Bottom: The experiment tracks the normalized amount of fluorescence in the bleach spot (orange dots) through time. Blue stars denote the amount of fluorescence in the circular bleaching spot corresponding to the above three panels. Fitting an appropriate PDE model of protein transport generates the black curve through the fluorescent recovery. }
\label{fig:FRAP_fit}
\end{figure*}

\textit{In vivo} imaging techniques can help quantify protein dynamics in cells. For example, microscopy techniques such as Fluorescence Recovery After Photobleaching (FRAP) and photoactivation provide time-series data quantifying the amount of mRNA at certain cellular locations throughout time in developing oocytes \cite{gagnon2013directional, powrie2016using}. FRAP has been particularly well-studied in conjunction with mathematical modeling. In FRAP experiments, fluorescent molecules in a small region of the cell are irreversibly bleached, and subsequent movement of the surrounding non-bleached fluorescent molecules into the photo bleached area is recorded, thus generating fluorescence intensity recovery curves as in Figure \ref{fig:FRAP_fit}.

Many methods have been developed to analyze FRAP data, typically using diffusion or reaction-diffusion equations and often making approximations based on whether diffusion or binding dominates the dynamics \cite{sprague2004analysis,ciocanel2017analysis}. In many biological applications, FRAP curves are fit to an exponential expression, which is only appropriate in scenarios where diffusion is very fast compared to the binding and timescale of the FRAP experiment (i.e., the reaction dominant scenario outlined in \cite{sprague2004analysis}). This is appropriate, for example, when quantifying RNA dynamics at the cell cortex, where RNA is likely to be a part of a highly stable and localized protein complex \cite{powrie2016using}. However, these approaches are unlikely to estimate parameters accurately in systems where active transport is a key mechanism, as well as in settings where simplifying assumptions do not apply. 

To connect intracellular transport models such as \eqref{eq:2-state} with microscopy data such as FRAP, we have previously proposed setting initial conditions $u(x,y,0)$ and $v(x,y,0)$ for transport PDE models based on FRAP postbleach intensity profiles from experiments \cite{ciocanel2017analysis}. Under reasonable assumptions about the experimental set-up, FRAP recovery curves then correspond to: 
\begin{align}
    \mathrm{FRAP}(t) = \int_{\mathrm{bleach \; spot}} (u+v)(x,y,t) dx dy \,.
\end{align}
Using numerical integration of the systems of advection-reaction-diffusion PDEs, we then used deterministic optimization methods to estimate the diffusivity $d$, active speed $c$, and transition rates $\beta_1$ and $\beta_2$ that characterize RNA dynamics \cite{ciocanel2017analysis} (see Figure~\ref{fig:FRAP_fit}). With limited knowledge about parameter values describing the dynamics in a living system such as the \textit{Xenopus} oocyte, and given that we are working with PDE systems, we found that parameter estimation can be challenging. In \cite{ciocanel2017analysis}, we first performed large parameter sweeps (which can nonetheless be computed in parallel) to determine good initial parameter guesses. A similar approach was used in \cite{sprague2004analysis}.

These parameter estimates can then be combined with analytically-derived effective transport quantities as in \eqref{eq:effective_qtities} and \eqref{eq:effective_qtities_2state} to provide insights for the biological problem. One advantage of the FRAP experimental technique is that it can be performed in different regions of a living cell, which enables direct observation of protein mobility throughout the oocyte. This has allowed us to combine experiments, parameter estimation, and analysis to suggest that mRNAs undergo faster uni-directional movement in the upper vegetal cytoplasm and higher spread or bidirectional transport in the lower vegetal cytoplasm of the oocyte, thus supporting previously-generated hypotheses from the lab \cite{ciocanel2017analysis}. 

Since several models may provide possible descriptions for a specific application, the question of model selection is also a highly relevant one in these applications. We have found it useful to determine how other quantities, such as the expected distances and times that the cargo spends on microtubules, depend on the model parameters. Evaluating these quantities at the estimated parameter values and comparing them with biologically-relevant time and length scales for microtubule travel can allow us to distinguish between and select appropriate modeling frameworks for the application \cite{ciocanel2017analysis}. 

On a personal note, this project also highlighted to me the importance of close collaborations with biologists. Continued conversations with lab members and an openness to model refinements are helpful in developing a common language with life sciences collaborators. In one specific instance, this helped me understand that aspects that may be simpler to model (such as assuming that the photobleach experiment is instantaneous, thus leading to simple initial conditions) are not appropriate and may even lead to misleading estimates and results.

\section*{Stochastic Modeling of Intracellular Transport}

An equivalent framework for analyzing the transport properties of intracellular cargo is to consider stochastic state-switching particle models. In these models, a particle (cargo protein) can switch between a finite number of states described by a continuous-time Markov chain $J(t)$, and its position is given by the real-valued process $X(t)$. For the simple model in Figure~\ref{fig:mRNA_oocyte_2state}B, $J(t)$ would only take the values $0$ for the diffusion state, and $1$ for the transport state. If we let $\vect{\pi}$ represent the stationary distribution of $J(t)$, this vector will consist of the long-term fractions of time spent in each state. The effective transport properties of cargo based on the dynamical systems approach in \cite{ciocanel2017analysis} and calculated in \eqref{eq:effective_qtities} can then be rewritten as:
\begin{align}\label{eq:effective_qtities_stoch}
    v_{\mathrm{eff}} &= \vect{c} \cdot \vect{\pi} \,, \nonumber \\
    \sigma_{\mathrm{eff}} &= \vect{d} \cdot \vect{\pi}  -  \vect{c} \cdot  \tilde{A}^{-1} (v \circ \vect{\pi} - v_{\mathrm{eff}} \vect{\pi}) \,,
\end{align}
where $\vect{c}$ is a vector of the active speeds in each state $\vect{c}=(c_0,...,c_n)^T$ and similarly $\vect{d}=(d_0,...,d_n)^T$ is a vector of the diffusivities in each state. Here $v \circ \vect{\pi}$ denotes the Hadamard product, or component-wise multiplication. 

This approach, as well as the quasi-steady state methods discussed in the previous section, rely on the Markovian structure of the dynamics, meaning that the state-switching process is a continuous-time Markov
chain with exponentially distributed durations in each state. The stochastic processes framework we extended in \cite{MR4153299} allows for more generalized random dynamics, which in particular need not assume that the particle spends an exponentially-distributed time in each state. 

In \cite{MR4153299}, we assume that the cargo switches between dynamic states at the random times $\{t_k : k \in \mathbb{N}\}$, and let $\{J_k : k \in \mathbb{N}\}$ denote the state during the $k$th time interval $ [t_{k-1}, t_k) $. We assume that the sequence $J_k$ is a time-homogeneous Markov chain and that it has a finite mean time of returning to each particular state (i.e., it is called positive recurrent). This assumption accurately describes many applications in cellular transport (though it does not describe settings with external stimuli). We therefore assume that the particle position dynamics have a regenerative structure, which repeatedly returns to a base state after some regeneration times. 

We can then answer questions about cargo transport properties in the context of renewal reward theory. Specifically, the dynamics can be considered as a sequence of independent cycles, characterized by returns to a chosen base state (such as the diffusive state). We denote the times to re-enter the base state as the regeneration times $T_n$ and define random variables for a generic cycle:
\begin{align}\label{eq:cycle_vars}
    \Delta T &= T_n - T_{n-1} \nonumber \,, \\
    \Delta X &= X(T_n) - X(T_{n-1})\,.
\end{align}
Based on the functional central limit theorem in \cite{MR2484222} as well as building on prior work on molecular motor systems \cite{MR2974479}, the effective velocity and diffusivity (i.e., the long-run properties of transport) are given by:
\begin{align}\label{eq:effective_qtities_renewal}
     v_{\mathrm{eff}} &= \frac{\mathbb{E}(X)}{\mathbb{E}(T)}  \nonumber \,, \\
    \sigma_{\mathrm{eff}} &= \frac{1}{2\mathbb{E}(\Delta T)} (\mathrm{Var}(\Delta X) + v_{\mathrm{eff}}^2\mathrm{Var}(\Delta T) -2v_{\mathrm{eff}} \mathrm{Cov}(\Delta X, \Delta T) ) \,.
\end{align}

This means that evaluating the long-run transport properties of cellular cargo reduces to understanding the moments of time and displacement random variables within a single generic cycle of the dynamics. To represent the contribution of events within a single cycle, we can re-write the cycle statistics as:
\begin{align}\label{eq:cycle_stats}
    \Delta T = \sum_{k=1}^\eta \tau_k \,, \;
    \Delta X = \sum_{k=1}^\eta \xi_k \,,
\end{align}
where $\eta$ is the number of events (steps) in a generic cycle. $\tau_k$ is the time spent in each step and $\xi_k$ is the corresponding spatial displacement of the particle in each step. We can think of these quantities as rewards collected in each step, which naturally depend on the biophysical state that the cargo is in at that time: $\tau_k = \sum_{j=0}^N \tau_k(j) 1_{\{J_k=j\}}$ and $\xi_k = \sum_{j=0}^N \xi_k(j) 1_{\{J_k=j\}}$. In the context of the 2-state model example in Figure~\ref{fig:mRNA_oocyte_2state}B, and assuming that the times spent by the cargo in each state are exponential, these times and displacements are given by:
\begin{align}
    \tau(0) &\sim \mathrm{Exp}(\beta_2)\,; \; \xi(0) \sim \sqrt{2d\tau(0)}Z \nonumber \,, \\
    \tau(1) &\sim \mathrm{Exp}(\beta_1)\,; \; \xi(1) \sim c \tau(1)\; \,.
\end{align}
As mentioned above, state $0$ corresponds to diffusion and state $1$ corresponds to active transport by motor proteins. The model parameters have the same meanings as in equations~\eqref{eq:2-state}, and $Z$ is an independent standard normal random variable. 

The analysis in \cite{MR4153299} proceeds by interpreting the time duration and spatial displacement in each state as rewards associated to that cargo state. An approach based on setting up moment-generating functions of the reward collected by the cargo until it returns to the base state allows to calculate the moments of the rewards needed to compute the effective transport properties in \eqref{eq:effective_qtities_renewal}. Interested readers can find additional details in \cite{MR4153299}. 

Applying this renewal reward framework to advection-reaction-diffusion models for the dynamics of mRNA protein cargoes yields effective transport quantities that agree with the results of dynamical systems analysis in the previous section (such as in \eqref{eq:effective_qtities_2state}). This is not surprising, since both approaches quantify the long-time behavior of the cargo and make the assumption that times spent in each state are exponential. The advantage of the stochastic processes approach is that times spent in each state can have non-exponential distributions, which may provide more accurate models in certain applications. We have found that, in models of \textit{in vitro} experiments studying how cargo is pulled by teams of molecular motors \cite{kunwar2010robust,muller2010bidirectional}, this method can complement costly numerical simulations and give insights into how effective transport quantities depend on individual motor properties.

\begin{figure}
\centering
\includegraphics[width=1\linewidth]{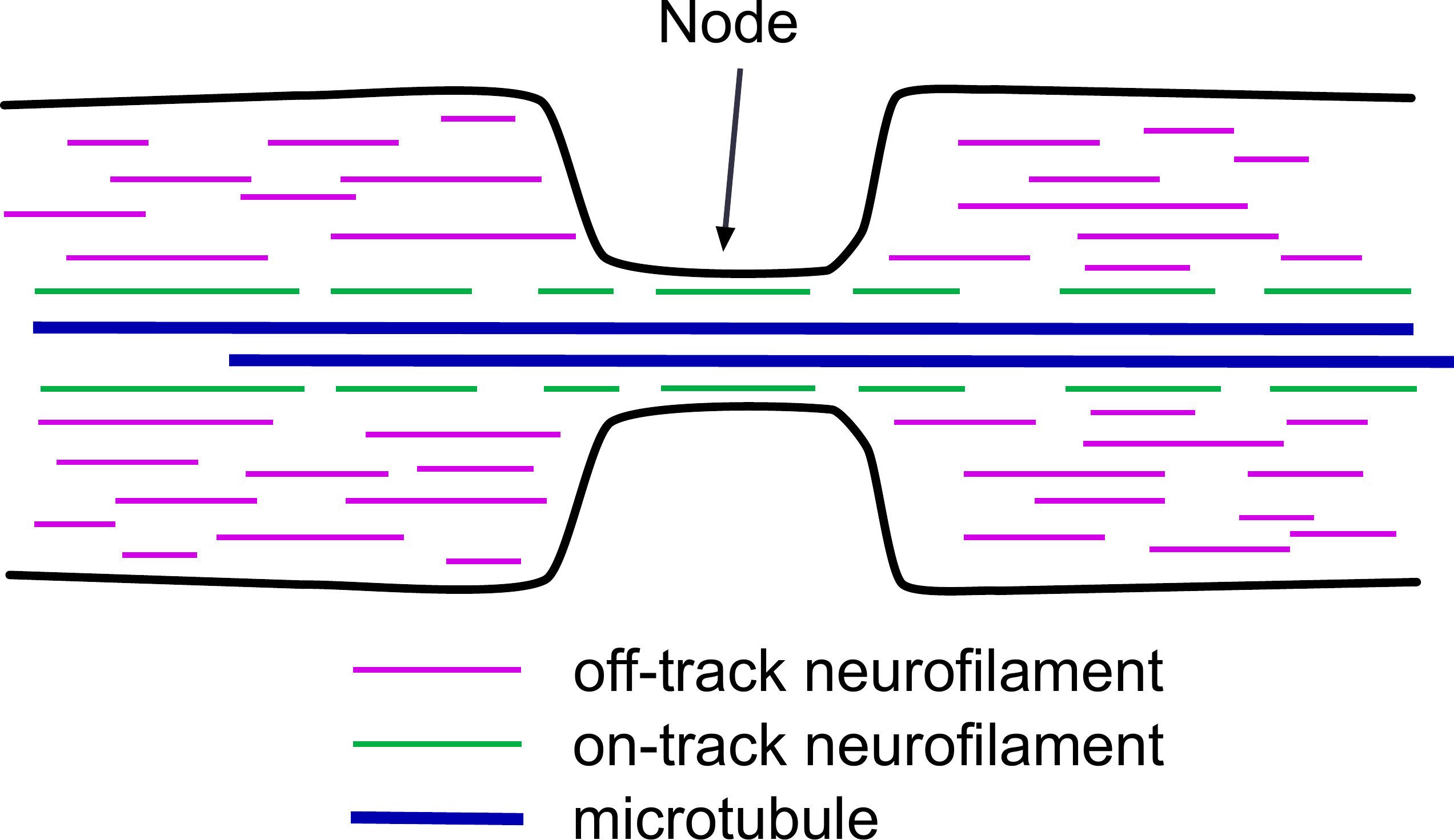}
\caption{Illustration of neurofilament transport along a neuronal axon. Microtubules (blue lines) go through axonal constrictions (nodes of Ranvier). Neurofilaments can switch stochastically between on-track states (green) and off-track states (purple). Experimental information about length distributions of neurofilaments can be incorporated into a six-state stochastic model and provide mechanistic insights on kinetic regulation of the dynamics at axonal nodes.}
\label{fig:NFs}
\end{figure}

In other applications, numerical simulations are most appropriate for investigating certain biological questions, especially if additional experimental data on protein sizes is available. One such example is illustrated by neurofilament transport, as in Figure~\ref{fig:NFs}. Neurofilaments are long and space-filling protein polymers that represent one of the most abundant cargoes in axons of neurons. Despite being an example of intermediate filaments, which are a component of the cytoskeleton of the cell, neurofilaments act as cargo and get transported along long microtubule filaments that ensure long-range transport along axons. Neurofilaments have been shown to move bidirectionally along microtubule tracks, illustrating so-called ``stop-and-go" dynamics characterized by short bouts of rapid movement on the filament track, interrupted by longer pauses off the microtubule tracks (see Figure~\ref{fig:NFs}). This dynamics has been accurately described by a stochastic model consisting of six kinetic states, each with different speeds and switching rates \cite{brown2005stochastic}.

Interestingly, in healthy neurons, neurofilaments successfully navigate axonal constrictions called nodes of Ranvier, which constitute potential bottlenecks for their transport. One mechanism that may explain this successful navigation and the observed acceleration of the polymers in nodes of Ranvier is that they may be physically closer to microtubule tracks in these nodes. Since the length distribution of neurofilaments has been recently characterized in model organisms such as mice, we complemented the existing stochastic model with length information \cite{ciocanel2020mechanism}, and used it to test this potential mechanism of transport and study how parameters describing the cargo dynamics may be regulated at the nodes (see Figure~\ref{fig:NFs}).

Specifically, we modeled each neurofilament as an individual linear structure, and carried out simulations of their stochastic cycling between the kinetic states and of their interaction with microtubules at node or non-node locations along a model axon. While more computationally costly than the stochastic process analysis previously described, this approach gave insight into the potential mechanistic behavior of the polymers at nodes and connected information about neurofilament numbers and lengths to axonal diameter. This framework also allowed us to validate our model with fluorescence pulse-escape experiments, where fluorescence of the protein is activated in a short region of the axon and tracked throughout time. This project is another example where carrying out mathematical model validation with data from experiments can support mechanistic understanding of cellular transport problems.

\section*{Cytoskeleton Geometry and Dynamics}

In the models described thus far, there is an implicit assumption that cargo proteins move along a single (one-dimensional) direction of active transport. While this may be a useful simplification when describing short-time experiments, such as microscopy experiments, it is no longer appropriate when considering long-term transport processes. For instance, mRNA localization takes 1-2 days in the developing frog oocyte, and its ultimate spatial localization, which impacts the fate of the developing organism, is likely to be significantly affected by the organization of microtubule filaments throughout the egg cell. In addition, cytoskeleton filaments have intrinsic dynamics of their own, which leads to interesting and functional cytoskeletal structures. The filament dynamics may be relevant to the transport of cargo, depending on the timescale of interest.

\subsection*{How about the Microtubule Transport Directions?}

Our analysis in the above sections assumes that transport is 1-dimensional (along a single filament track) and that diffusion is 2-dimensional. However, microtubules have complex orientation and organization in cells. In frog oocytes, they have a radial orientation bias, while in neurons, microtubules are mostly oriented parallel to each other. To allow for explicit analysis, we have considered a scenario of parallel microtubule tracks, which models filaments in axons or dendrites of neurons. 

Assuming that the density of filaments is allowed to vary with spatial dimension $x$ (assumed to be non-dimensionalized to $x \in [0,1]$), we and others have shown that this scenario is equivalent to the advection-reaction-diffusion PDE system \eqref{eq:gen}. The difference is that the reaction-rate matrix $A = A(x)$ now depends on space and accounts for the availability of microtubules at each location \cite{MR3892411,MR3328155}. Using Lyapunov-Schmidt reduction theory with some additional requirements on the resulting operators, we found that the asymptotic solution at large time is similar to \eqref{eq:ul}:
\begin{align}\label{eq:u_mts}
   \vect{u}(x,y,t) = \frac{1}{\sqrt{2\pi \sigma_{\mathrm{eff}} t}} e^{-\frac{(y+v_{\mathrm{eff}}t)^2}{2\sigma_{\mathrm{eff}} t}} \vect{u_0}(x)\,.
\end{align}
This analysis assumed two-dimensional diffusion \cite{MR3892411}. While the expression for the effective velocity $v_{\mathrm{eff}}$ is identical to the expression in \eqref{eq:effective_qtities}, the effective diffusivity now takes the form \cite{MR3892411}:
\begin{align}\label{eq:effective_diff_mts}
    \sigma_{\mathrm{eff}} &= \frac{\langle \vect{\psi_0},D\vect{u_0}\rangle}{\langle \vect{\psi_0}, \vect{u_0} \rangle} + \frac{\langle \vect{\psi_0},C\vect{w_0}(x)\rangle}{\langle \vect{\psi_0}, \vect{u_0} \rangle}\,.
\end{align}
Here $\langle f(x),g(x) \rangle = \int_0^1 f(x)g(x) dx$ and $\vect{w_0}(x)$ can be found by solving
\begin{align}\label{eq:w0_mts}
    (D\partial_x^2 + A(x)) \vect{w_0}(x) + C \vect{u_0}(x) - \frac{\langle \vect{\psi_0},C\vect{u_0}\rangle}{\langle \vect{\psi_0}, \vect{u_0} \rangle} \vect{u_0}(x) = 0 \,.
\end{align}
When applying this approach to the simple 2-state model of cargo switching between active transport and diffusion (Figure~\ref{fig:mRNA_oocyte_2state}B), we found that incorporating the spatial dependence of the microtubule structure enhances the effective diffusivity, and thus provides a more accurate prediction of the spread of the cargo at large time \cite{MR3892411}. 

Through this approach, we can fully derive the effective velocity and diffusion of the cargo under no assumptions on the density
of the microtubules, the magnitude of the reaction rates, or of the diffusion coefficient characterizing the particle dynamics. Transport in the context of more complex geometries has been studied using quasi-steady-state analysis and under assumptions of fast switching rates, slow diffusion, and sufficiently small filament density in \cite{MR3328155}, for applications including filaments growing from a microtubule aster in neuron growth cones, or filaments that nucleate from membrane sites in budding yeast. On the other hand, incorporating the dependence on spatial microtubule tracks in the stochastic regenerative cycle framework would require additional assumptions on the spatial scale of transport in a regeneration cycle \cite{MR4153299}.

\subsection*{Road networks have dynamics of their own}

In addition to forming complex geometries in some cells, cytoskeleton filaments are also dynamic. Microtubules have a dynamic plus end where they grow and shrink, while the opposite end is called the minus end. Similarly, actin filaments have polymerization and depolymerization dynamics at their plus and minus ends. 

\begin{figure}
\centering
\includegraphics[width=0.7\linewidth]{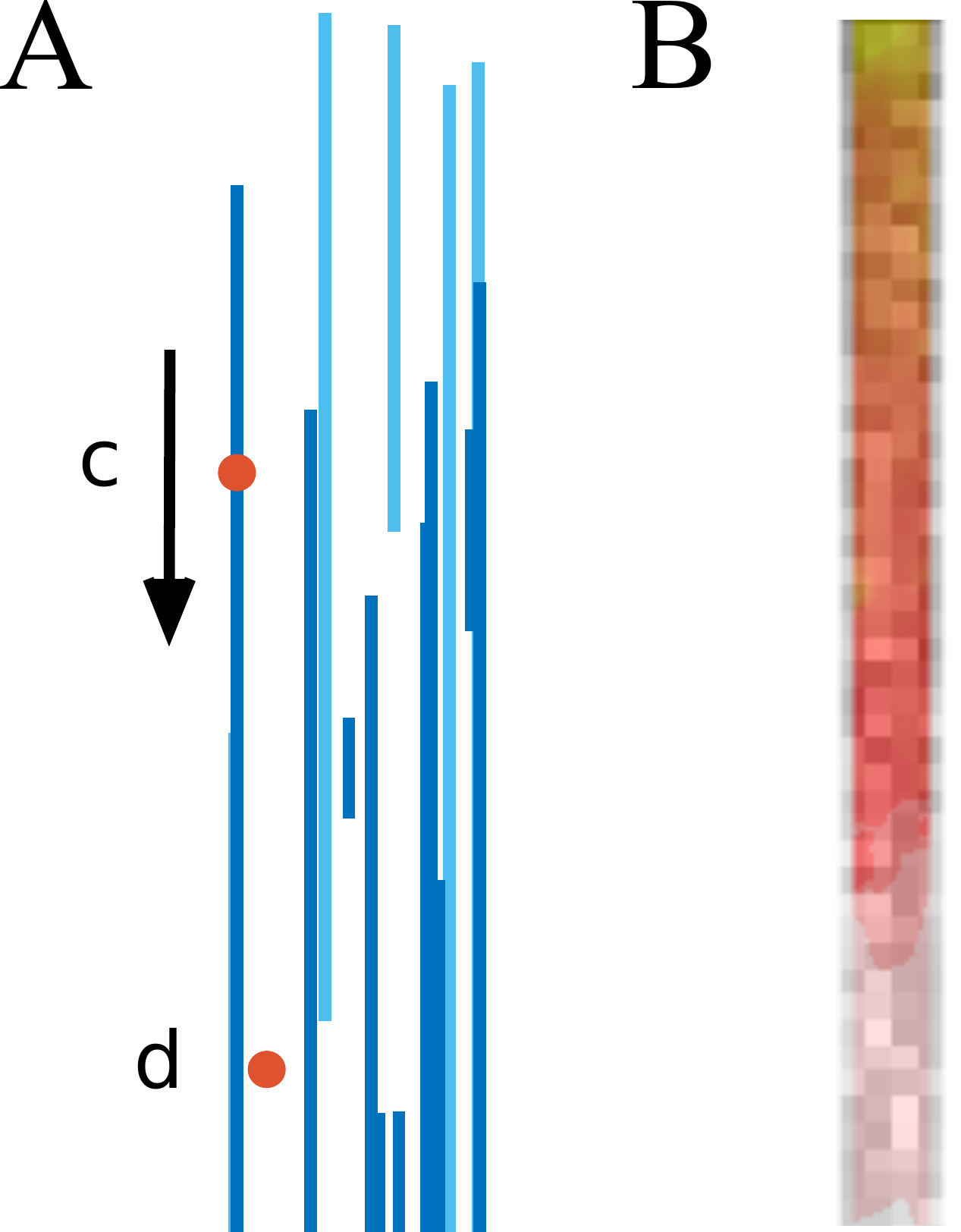}
\caption{(A) Sample model of parallel microtubules along a section of a neuronal dendrite. Dark blue lines denote microtubules with their plus ends oriented down, and light blue lines correspond to microtubules with the opposite orientation. Cargo (red particles) may switch between transport along microtubule tracks and diffusion off them.  (B) Numerical simulation of protein cargo localization predicted by the 2-state model in Figure~\ref{fig:mRNA_oocyte_2state}(B) with advection restricted to microtubule structures as in (A). }
\label{fig:localization}
\end{figure}

Dynamic microtubules undergo assembly and disassembly under various conditions in cells, so that the whole microtubule network may change within minutes, as observed in fruit fly oocytes \cite{trong2015cortical}. Extending the numerical simulation approach in \cite{trong2015cortical}, we have previously modeled random microtubule networks with a radial bias, as informed by experimental observations in frog oocytes \cite{MR3892411}. We constructed several such microtubule structures to account for the changes in the cytoskeletal network in time, and numerically simulated PDE models of protein transport on these networks using finite-volume discretization in 2-dimensional space. While this setting is challenging to study using analytical techniques, these numerical simulations allowed us to use parameters estimated from microscopy experiments and to generate hypotheses about anchoring of mRNA at the oocyte cortex \cite{MR3892411}, which inform future biological experiments. This approach may also prove useful in understanding how different orientations of microtubules, as in Figure~\ref{fig:localization}(A), impact transport of protein cargoes in neurons, as predicted by numerical simulations (Figure~\ref{fig:localization}(B)).

It is clear that dynamic cytoskeleton networks play an important role in directing protein transport in cells. Additionally, the dynamic nature of cytoskeletal filaments can also result in cell shape changes or regulate the mechanical profile of the cell. For example, actin filaments are highly dynamic and can quickly respond to internal and external cues in the cell cortex. Understanding filament reorganization has led to other research directions, such as developing stochastic agent-based models incorporating chemical and mechanical protein interactions \cite{popov2016medyan}. The complex spatiotemporal data emerging from this work has also motivated the development of novel data analysis techniques drawing on network theory, spatial statistics, or topological data analysis \cite{popov2016medyan,MR4200885}.

\section*{Outlook and Open Challenges}

I have become interested in mathematical cell biology while working closely with life sciences labs throughout my training. This is an exciting field to contribute to, given that our understanding of biological processes is constantly getting updated, both through innovation in experiments and through model refinements and advances in mathematical analysis and computational techniques. A key feature of working in this field is the importance of being open to changing modeling frameworks as additional biological information becomes available.

As an example, when I started working on mathematical modeling of mRNA transport, we thought of mRNA cargo as molecules that get transported along microtubule filaments in developing oocytes. In the past years, advances in imaging and microscopy have helped our collaborators understand that mRNAs actually organize in granules, as observed in other cells \cite{neil2021bodies}. This opens up interesting questions about the within-granule dynamics and organization, which is characterized by phase separation and has been recently modeled in \cite{MR4181743}. Understanding effective transport of such granules will likely introduce new mathematical and computational challenges. New methods are also being developed for visualizing cytoskeleton networks, which will ultimately improve our understanding of cytoskeleton dynamics, as well as open up new directions for data-driven mathematical modeling. Connecting PDE and stochastic models with experimental data remains a significant challenge in this field and drives questions about parameter estimation and identifiability given existing modeling frameworks. Close collaborations between applied mathematicians and experimentalists will be key to addressing protein transport questions that influence our understanding of healthy cellular function and organization.

\bibliography{references}

\end{document}